\documentclass[preprint,12pt]{elsarticle}

\usepackage{amssymb}
\usepackage{amsthm}
\usepackage{amsmath}
\usepackage{mathrsfs}
\usepackage{tikz}
\usetikzlibrary{calc,intersections,through,backgrounds}
\usepackage{graphicx}
\usepackage{caption}
\usepackage{subcaption}

\newtheorem{thm}{Theorem}
\newtheorem{lem}[thm]{Lemma}
\newdefinition{defn}{Definition}
\newproof{pf}{Proof}
\journal{Journal of Discrete Mathematics}

\begin{document}

\begin{frontmatter}


\title{The Vapnik-Chervonenkis dimension of cubes in $\mathbb R^d$}
\author{Christian J. J. Despres\tnoteref{presentaddress}}
\tnotetext[presentaddress]{Present Address: University of
    Toronto, 40 St. George Street, Toronto, ON, M5S 2E4}
\ead{christian.despres@mail.utoronto.ca}
\address{University of Ottawa, 585 King Edward, Ottawa, ON, K1N 6N5}

% \title{}


% \author{}

% \address{}

\begin{abstract}
  The Vapnik-Chervonenkis (VC) dimension of a collection of subsets of a set is
  an important combinatorial concept in settings such as discrete geometry and
  machine learning. In this paper we prove that the VC dimension of the family
  of $d$-dimensional cubes in $\mathbb R^d$ is $\lfloor(3d+1)/2\rfloor$.
\end{abstract}

\begin{keyword}


\MSC[2010] 03E05 \sep Vapnik-Chervonenkis dimension \sep discrete geometry 
\end{keyword}

\end{frontmatter}
\section{Introduction}
\label{intro}
The Vapnik-Chervonenkis (VC) dimension of a collection $\mathscr C$ of subsets
of a space $\Omega$ is one measure of the complexity of $\mathscr C$. Introduced
in \cite{vapnik-chervonenkis}, it has found wide application in areas such as
machine learning, where it is used to gauge the capacity of a model to represent
sample data (see, e.g., \cite{vapnik} and \cite{devroye-lugosi}, or
\cite{floyd-warmuth} for work on the related problem of sample compression). The
VC dimension of many natural collections of subsets of Euclidean space has been
determined. For instance, it is a standard result that the VC dimension of
balls, half-spaces, and boxes in $\mathbb R^d$ is $d+1$, $d+1$, and $2d$,
respectively. However, the VC dimension of cubes, or more generally the balls
according to the $\ell^p$ norm on $\mathbb R^d$ for $p \neq 2$, has not
previously been calculated. In this paper we will show that the VC dimension of
cubes (the balls of the $\ell^\infty$ norm) in $\mathbb R^d$ is
$\lfloor(3d+1)/2\rfloor$. A remaining question is the VC dimension of the balls
of other $\ell^p$ norms. While the VC dimension of such collections of balls for
$p$ sufficiently large should be at least that of cubes, the precise values of
the VC dimension and its behaviour as $p$ changes is not known.

To define the VC dimension, let $\Omega$ be a set and $\mathscr C$ a collection of subsets of $\Omega$. If $S$ is a subset of $\Omega$ and $A$ a subset of $S$, then we say that $\mathscr C$ \emph{carves out} $A$ from $S$ if there is a set $C$ in $\mathscr C$ such that $S \cap C = A$. The set $S$ is \emph{shattered} by $\mathscr C$ if every subset of $S$ is carved out by $\mathscr C$. With these definitions, we have the following.

\begin{defn}
The VC dimension of a collection $\mathscr C$ of subsets of $\Omega$ is the supremum of the cardinalities of finite sets that are shattered by $\mathscr C$.
\end{defn}
To prove that the VC dimension of cubes in $\mathbb R^d$ is
$\lfloor(3d+1)/2\rfloor$, we first show that no set of larger cardinality is
shattered by cubes, and then we construct sets in $\mathbb R^d$ of the
appropriate size that are shattered by cubes. Before we proceed, we remark that
we will refer to the $i$th coordinate of a point in $\mathbb R^d$ with a
superscript.

\section{Establishing the upper bound on the VC dimension}
We now prove that the VC dimension of cubes is at most $\lfloor(3d+1)/2\rfloor$.

\begin{thm}
  Any set in $\mathbb R^d$ that is shattered by cubes has at most $\lfloor(3d+1)/2\rfloor$ points.
\end{thm}
\begin{proof}
Let $S$ be a subset of $\mathbb R^d$ with $d+n$ points, and suppose that $S$ is
shattered by cubes. In order to establish the upper bound, for every axis $i$ we
pick points $l_i$ and $u_i$ in $S$ whose $i$th coordinates are minimal and
maximal, respectively, among the $i$th coordinates of points in $S$. Observe
that every point in $S$ must appear somewhere in the list $(l_1, u_1), \dotsc,
(l_d, u_d)$ of pairs of extremal points. For if that were not the case, and a
point $a$ in $S$ did not appear, then the subset $S \setminus \{a\}$ could not
be carved out by a cube.

Now, let $k$ be the number of points in $S$ that appear precisely once in our
list. These $k$ points are distributed over $2d$ positions in the list of
extrema, leaving $2d-k$ positions in the list to be filled by the $d+n-k$ points
that appear at least twice. So $2d-k \geqslant 2(d+n-k)$, implying that $k
\geqslant 2n$.

Finally, assume towards a contradiction that $k \geqslant d+2$. By the
pigeonhole principle there must be distinct axes $i$ and $j$ such that the
points $l_i$, $u_i$, $l_j$, and $u_j$ appear precisely once in the list of
extrema. Since we can carve out $S\setminus\{l_i, u_i\}$ with a cube, we can
find two closed intervals of the same length such that the first contains the
$j$th coordinate of each point in $S$ and the second contains the $i$th
coordinate of each point in $S$ except for those of $l_i$ and $u_i$. Thus
$|u_i^i - l_i^i| > |u_j^j - l_j^j|$. By repeating this argument with the roles
of $i$ and $j$ interchanged we obtain the contradictory inequality $|u_j^j -
l_j^j| > |u_i^i - l_i^i|$. Thus $d+1 \geqslant k \geqslant 2n$, giving the upper bound.
\end{proof}

\section{Establishing the lower bound on the VC dimension}
In order to construct sets in $\mathbb R^d$ of size $\lfloor(3d+1)/2\rfloor$ that can be shattered
by cubes, we need some preliminary definitions. Consider the collection of sets
$\mathscr I$ consisting of $\mathbb R$, all closed left-infinite intervals $(-\infty,a]$ with $a \geqslant 0$,
and right-infinite intervals $[a,\infty)$ with $a \leqslant 0$. Call a product
$I = I^1 \times \dotsb \times I^d$ of intervals in $\mathscr I$ a \emph{constraint} in $\mathbb R^d$. We say that $I$ is
\emph{right-exclusive} if $I^d$ is half-infinite to the left (that is, of the
form $(-\infty,a]$), \emph{left-exclusive} if $I^d$ is half-infinite to the right, and \emph{inclusive} if $I^d$ is $\mathbb R$. Observe that given a constraint $I$ and a bounded set $S$, it is possible to carve out $S \cap I$ with arbitrarily large cubes. So if we can shatter a set $S$ with constraints, then we may do so with cubes. 

To construct shattered sets of the appropriate size, we will need to construct
sets in each dimension $d$ that are slightly smaller than required but can be
shattered by constraints in a particularly ``nice'' way. The construction of
these sets will be recursive in $d$ (see Lemmas \ref{addtwo} and \ref{addone}), and the niceness property that is captured
in Definition \ref{niceness} will allow the recursion to continue.

\begin{defn}
\label{niceness}
Let $S$ be a subset of $\mathbb R^d$. Say that $S$ is \emph{accessible} if every subset of $S$
may be carved out by an inclusive constraint. Say that $S$ is \emph{weakly accessible} if
every subset of $S$ may be carved out by two constraints, one of which is
left-exclusive and the other of which is right-exclusive.
\end{defn}

\begin{lem}
\label{addtwo}
If $S$ is an accessible subset of $\mathbb R^d$ then we may adjoin two points to $S$ to
obtain a set $S'$ of points shattered by cubes. We may also embed $S'$ in
$\mathbb R^{d+1}$ to obtain a weakly accessible subset of $\mathbb R^{d+1}$.
\end{lem}

\begin{proof}
  Given such an $S$, let $x$ and $y$ be points in $\mathbb R^d$ with $x^d$
sufficiently negative and $y^d$ sufficiently positive, so in particular that $x$
and $y$ are minimal and maximal, respectively, on the $d$th axis, and with
$x^j=y^j=0$ for $j < d$. Let $S' = S \cup \{x,y\}$. To verify that $S'$ is
shattered by cubes, let $A$ be a subset of $S$ and $I$ an inclusive constraint
carving out $A$. Then $I$ carves out $A \cup \{x,y\}$ from $S'$. We may carve
out the subsets $A \cup \{x\}$ and $A \cup \{y\}$ by taking the constraint $I$
and replacing $I^d$ with an appropriate half-infinite interval from $\mathscr
I$. Finally, we may carve $A$ from $S'$ by carving $A$ from $S$ with a
small cube.

  Now embed $S'$ in $\mathbb R^d$ by taking each point in $S'$ and duplicating
its $d$th coordinate to serve as its $(d+1)$st coordinate. The resulting set is
weakly accessible; the additional axis $d+1$ allows us to include and exclude
the images of $x$ and $y$ from any subset we wish to carve out, and to do so
with both left- and right-exclusive constraints.
\end{proof}

\begin{lem}
\label{addone}
If $S$ is a weakly accessible subset of $\mathbb R^d$ then we may adjoin one
point to $S$ to obtain a set $S'$ that is shattered by constraints, and hence by
cubes. We may also embed $S'$ in $\mathbb R^{d+1}$ to obtain an accessible
subset of $\mathbb R^{d+1}$.
\end{lem}

\begin{proof}
Let $x$ be a point in $\mathbb R^d$ whose $d$th coordinate is strictly smaller than
those of the points in $S$, and whose $j$th coordinate is $0$ for each $j < d$.
Define $S' = S \cup \{x\}$, and consider a subset $A$ of $S$. A left-exclusive
constraint $I$ that carves out $A$ from $S$ (with the left endpoint of $I^d$ not
too small) will also carve out $A$ from $S'$, while a right-exclusive constraint
that carves out $A$ from $S$ will carve out $A \cup \{x\}$ from $S'$. So $S'$ is
shattered by constraints.

As in the previous lemma, embed $S'$ into $\mathbb R^{d+1}$ by duplicating the $d$th
coordinate of each point in $S'$. Since $S'$ is shattered by constraints, the
resulting set in $\mathbb R^{d+1}$ is accessible; we do not need to use the
$(d+1)$st axis to exclude any points when carving out subsets from the image of $S'$.
\end{proof}

Note that the empty set in $\mathbb R^1$ is accessible. By alternately applying
Lemmas \ref{addtwo} and \ref{addone} to this base case we obtain Theorem
\ref{construction}.

\begin{thm}
\label{construction}
For each $d$ there is a set $S$ in $\mathbb R^d$ with $\lfloor(3d+1)/2\rfloor$ points that is
shattered by cubes.
\end{thm}

\bibliographystyle{elsarticle-num} 
\bibliography{myrefs}

\end{document}